# A Review on Constraint Handling Techniques for Population-based Algorithms: from single-objective to multi-objective optimization

Iman Rahimi[1,] **Amir H. Gandomi**[1*,] **Fang Chen**[1,] **and Efrén Mezura-Montes**[2]

[1] Data Science Institute, *Faculty of Engineering & Information Technology, University of Technology Sydney, Australia*; iman83@gmail.com , Gandomi@uts.edu.au, fang.chen@uts.edu.au

[2] *Artificial Intelligence Research Institute, University of Veracruz, MEXICO*, emezura@uv.mx

**\*** Correspondence: Gandomi@uts.edu.au;

**Abstract-** Most real-world problems involve some type of optimization that is often constrained. Numerous researches have investigated several techniques to deal with constrained single-objective and multi-objective evolutionary optimization in many fields, including theory and application. Also, this presented study provides a novel analysis of scholarly literature on constraint handling techniques for single-objective and multi-objective population-based algorithms according to the most relevant journals, keywords, authors, and articles. The paper reviews the main ideas of the most state-of-the-art constraint handling techniques in multiobjective population-based optimization, and then the study addresses the bibliometric analysis in the field. The extracted papers include research articles, reviews, book/book chapters, and conference papers published between 2000 and 2020 for the analysis. The results indicate that the constraint handling techniques for multiobjective optimization have received much less attention compared with single-objective optimization. The most promising algorithms for such optimization were determined to be genetic algorithms, differential evolutionary algorithms, and particle swarm intelligence. Additionally, "Engineering," "Computer Science," and " Mathematics" were identified as the top three research fields, in which future research work is anticipated to increase.

## 1. Introduction

Most real-world problems are considered as multi-objective optimization problems (MOOPs). No single solution exists for an MOOP, instead different solutions generate trade-offs for different objectives. Furthermore, MOOPs arise in a natural fashion in most disciplines, and solving them has been a challenging issue for researchers. Although a variety of methods have been developed in Operations Research and other fields to address these problems, there is an urgent need for alternative approaches because of the complexities of their solution [1]–[3]. Evolutionary computation (EC) methods have been identified as more effective methods to handle this limitation and are suitable for solving MOOPs, for which the form of the Pareto-optimal front (discontinuity, nonconvexity, etc.) is not important [4], [5]. Moreover, most multiobjective evolutionary algorithms (MOEAs) use the dominance concept [6], [7].

To solve the constrained optimization of all real-world problems, constrained evolutionary algorithm optimization (CEAO) implements an evolutionary algorithm (EA) combining with a constraint handling technique (CHT). In a work by [8], an infeasible individual will be divided into different categories based on their distances to the feasible region, and ranking will be conducted according to the classes. [9] introduced an approach that assigns high and low priorities to constraints and objective functions, respectively. The authors of [10] proposed a CHT that only considers the inequality constraints, wherein the algorithm uses tournament selection that has better convergence properties in comparison to the proportionate selection operator [11]. However, the latter algorithm employs niche count for all populations that may increase the complexity of the computation.

The authors of [12] Introduced a novel approach that ignores any solution that violates any of the assigned constraints. [13] First proposed the use of a genetic algorithm (GA) population-based approach plus a controlled mutation operator to keep diversity among feasible solutions. The work of [14] proposed a CHT where three different non-dominated rankings of the population are performed, using objective function values, different constraints, and a combination of all objective functions and values. Although the technique can handle the infeasible solutions carefully and maintain diversity in the population, the algorithm performs poorly in terms of choosing parameter values and is computationally expensive. The authors of [15] developed an EA based on the nondominated sorting concept that uses the min-max formulation for constraint handling.

The authors of [16] Ran the simulation of the non-dominated sorting genetic algorithm II (NSGAII) on a seven-constrained nonlinear problem, which exhibited better performance than Ray-Tai-Seow's algorithm. The study of [17] conducted an overview and analysis of the most popular CHTs using EAs along with pros and cons. The authors of [18] Combined a penalty function and multiobjective optimization technique, in which the ranking scheme is borrowed from latter technique. The authors of [19] suggested two approaches, namely Objective Exchange Genetic Algorithm and Objective Switching Genetic Algorithm, for solving constrained MOOPs. A new partial order relation from the constraint MOOPs was proposed by [20], under which the Pareto optimum set satisfies the constraints. The research of [21] introduced the Blended Space EA, which uses a rank obtained by blending an individual's rank in the objective space to check dominance.

The authors of [22] Introduced a two-phase algorithm, which separates the objective function and constraints. The work of [23] introduced a MOO-based EA (Cai and Wang method), abbreviated as CW, in addition to three other models for constrained optimization. In the proposed approach, the simplex crossover was used to enrich the exploitation and exploration abilities. The authors of [24] proposed an EA based on evolutionary strategy for constrained multiobjective optimization problems. The method uses a min-max formulation for constraint handling in which feasible individuals and infeasible individuals evolve toward Pareto optimality and feasibility, respectively. The study of [25] suggested Pareto Descent Repair (PDR) to search for feasible solutions out of infeasible individuals. The authors of [26] proposed an adaptive tradeoff model for constrained evolutionary optimization to address three main issues: evaluating an infeasible solution in case the population contains only infeasible solutions; achieving a balance between feasible and infeasible individuals when the population contains both solutions; and selecting the feasible solution in case the population possesses only feasible solutions.

The authors of [27] Suggested a heuristic hybrid of particle swarm optimization (PSO) and ant colony optimization for the optimum design of trusses, which showed to handle the problem-specific constraints using a fly-back mechanism. The work of [28] suggested an infeasibility-driven EA (IDEA), which is able retain a proportion of infeasible solutions among the population members and preserve diversity compared to NSGAII. The authors of [29] investigated an EA solution for approximate Karush-Kuhn-Tucker (KKT) conditions of smooth problems. The results on some test problems indicate that EA's operators lead the search process to a point close to the KKT point. The authors of [30] discussed the most important techniques, many of which were previously proposed [17], [31]. The previous work also addressed some state-of-the-art constrained handling techniques, including feasibility rules based on GA [13], epsilon constrained method [32], penalty functions [33], [34], and ensemble of constraint-handling methods [35], [36].[37] Introduced an evolutionary scheme for handling boundary constraints and combined it with differential evolution (DE) and compared the proposed method with other boundary constraints handling techniques. The results indicated the proposed approach is much better than the existing methods.

The authors of [38] Developed a water cycle algorithm, inspired by observations of the water cycle process that could be applied to a number of constraint optimization problems. [39] Introduced a population-based algorithm based on the mine blast explosion concept then applied the proposed approach to some constraint optimization problems in comparison to other well-known optimizers. The authors of [40] used a constraint consensus method that helps an infeasible individual to move towards the feasible region and then combined the method with a memetic algorithm. The research conducted by [41] developed a feasible-guiding strategy to guide the evolution of individuals, in which a revised objective function technique with feasible guiding strategy based on NSGA-II is introduced to handle constrained MOOPs. The study proposed by [42] proposed a class of constraint handling strategies, which infeasible individuals are repaired when they are considered in the search space and explicitly preserve feasibility of the solutions.

The authors of [43] Used a hybrid of PSO and GA to improve the balance between exploration and exploitation by using genetic operators, namely crossover and mutation in PSO. A few years later, [44] extended the parameter-less CHT so as to provide a balance between the feasible and infeasible solutions in a GA population. The authors of [45] proposed a new approach, known as boundary update (BU) technique, which is able to handle constraints directly (i.e. updating variable bounds) and tested the proposed method on several constrained optimization problem. The method proposed by [45] possesses the potential to couple with MOEA.

It is noteworthy to mention that the majority of the mentioned studies focused on CHTs for single objective optimization with little attention towards multiobjective optimization. This is attributed to fact that most constraint handling methods developed for single objective optimization could be modified for multiobjective optimization as well [30].

To attain a better understanding of the research field and to provide new insights from relevant publications, this work aimed to answer the following questions:

- RQ1: What are CHTs, and how are they important?
- RQ2: What are the disadvantages of the different CHTs?
- RQ3: What are the main topics and keywords regarding constraint handling techniques?
- RQ4: Which journals have the most contributions in the field? Who are the best researchers in the area, and what countries are they from?
- RQ5: What are the basic statistics of constraint handling techniques for multiobjective population-bassed algorithms?
- RQ6: What are the most active countries and affiliations in the field?
- RQ7: What are the current gaps and future trajectory in the area?

The reminder of the study is as follows. Section 2 describes the research methodology. Section 3 presents the CHTs in EAs. Section 4 describes the CHTs in nature-inspired algorithms. Section 5 addresses the CHTs in multiobjective genetic algorithms. Section 6 presents a summary of novel approaches between 2020 and 2021. Section 7 discusses the scientometric analysis. Section 8 provides a summary of the study along with recommendations for future research. Concluding remarks are offered in the last section.

## 2. Research Methodology

The research procedure in this work was divided into five stages (Figure 1). In the first stage, documents from databases were gathered from Scopus and Web of Science (WOS). For this aim, the authors used special keywords, namely (TITLE-ABS-KEY (constrained AND multi AND objective AND evolutionary AND optimization) OR TITLE-ABS-KEY (constraint AND handling AND multi AND objective AND evolutionary AND optimization) OR TITLE-ABS-KEY (constrained AND multi AND objective AND swarm AND optimization) OR TITLE-ABS-KEY (constraint AND handling AND multi AND objective AND swarm AND optimization) to find the related articles published as of May 4, 2021. Supplementary A and B present the data extracted from Scopus and WOS, respectively. Since some of the articles were duplicates, they were identified and removed from the library in stage 2 using Mendeley as a powerful reference manager. Also, in stage 2, some research questions for this study were designed. An overview was initiated in stage 3 with a general illustration of the basic concepts of CHTs and comparison of methods. In stage 4, a social network analysis was performed to provide a scientometric analysis of the documents using VOSviewer [46], [47] and RStudio, which have been identified as powerful tools for scientometric analysis. Also, some interesting analytical features, such as number of pages and authors per article, were conducted in this stage. Moreover, stage 4 required several steps, including co-occurrence, co-authorship, citation, and citation network analyses. In the last stage, the results were obtained to formulate a discussion to answer the proposed research questions. Stage 5 required preparing the findings, identifying important gaps, and determining future research directions.

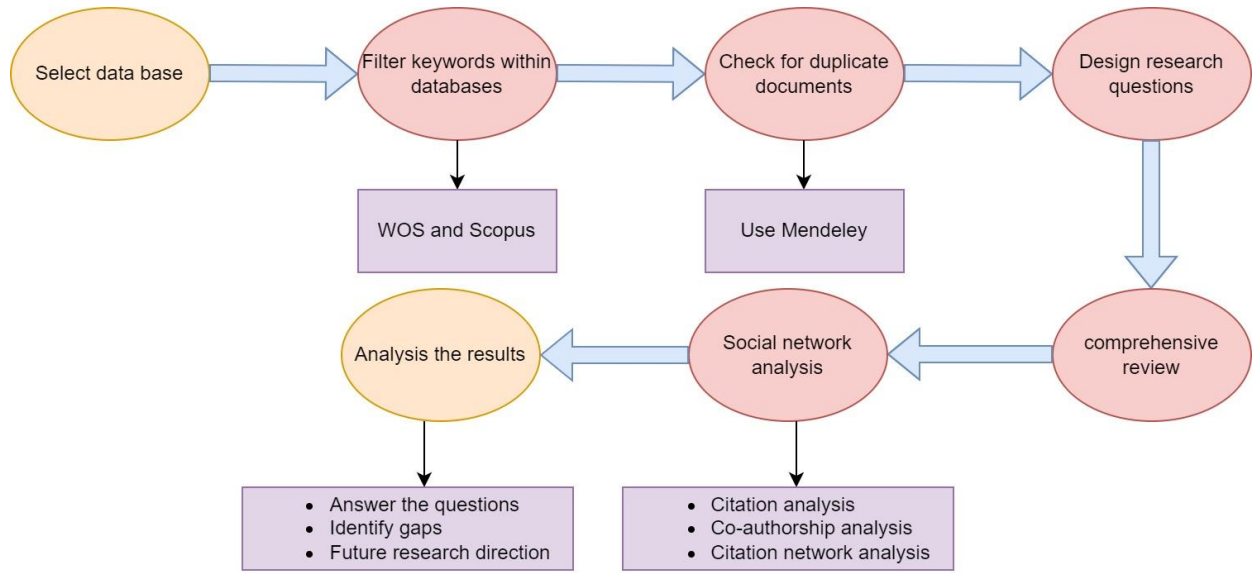

*Figure 1 Research Procedure*

## 3. Constraint Handling Methods in Evolutionary Algorithms (RQ1)

Almost all real-world problems are considered as constraint problems. A general form of a constrained multi-objective optimization problem(CMOOP) is described as follows (Equations 1-3):

Maximize (Minimize)　　　　　　　　　　(1)
F(x)=($f_1(x), ..., f_t(x)$)
s.t.
$h_i(x) = 0 \quad i = 1, ...., n$　　　　　　(2)
$g_j(x) \leq 0 \quad j = 1, ..., m$　　　　　　(3)

Where F(x) is the objective vector; and t, n and m are the number of objective function, equality and inequality constraints, respectively. There is no single solution for a MOOP that simultaneously optimizes each objective, instead there exists a number of Pareto optimal solutions. A Pareto front of possible solutions is called optimal or nondominated if improving anyone's objective further would lead to a decrease in other objectives. According to previous surveys [17], [31], a simple taxonomy of the constraint handling methods in nature-inspired optimization algorithms is as follows:

1- Penalty functions methods
2- Decoders
3- Special operators
4. Separation techniques

The first and fourth techniques are discussed in details later in the paper. As an example of decoders,[48] proposed a homomorphous mapping (HM) method between an n-dimensional cube and feasible space. The feasible region can be mapped onto a sample space where a population-based algorithm could run a comparative performance [48]–[51]. However, this method requires high computational costs. A special operator is used to preserve the feasibility of a solution or move within a special region [52]–[54]. Nevertheless, this method is hindered by the initialization of feasible solutions in the initial population, which is challenging with highly-constrained optimization problems. In addition, the authors of [55] presented a taxonomy of CHTs in MOEA as follows (Figure 2):

- Penalty functions
- Separation method
- Retaining the infeasible solutions
- Hybrid methods

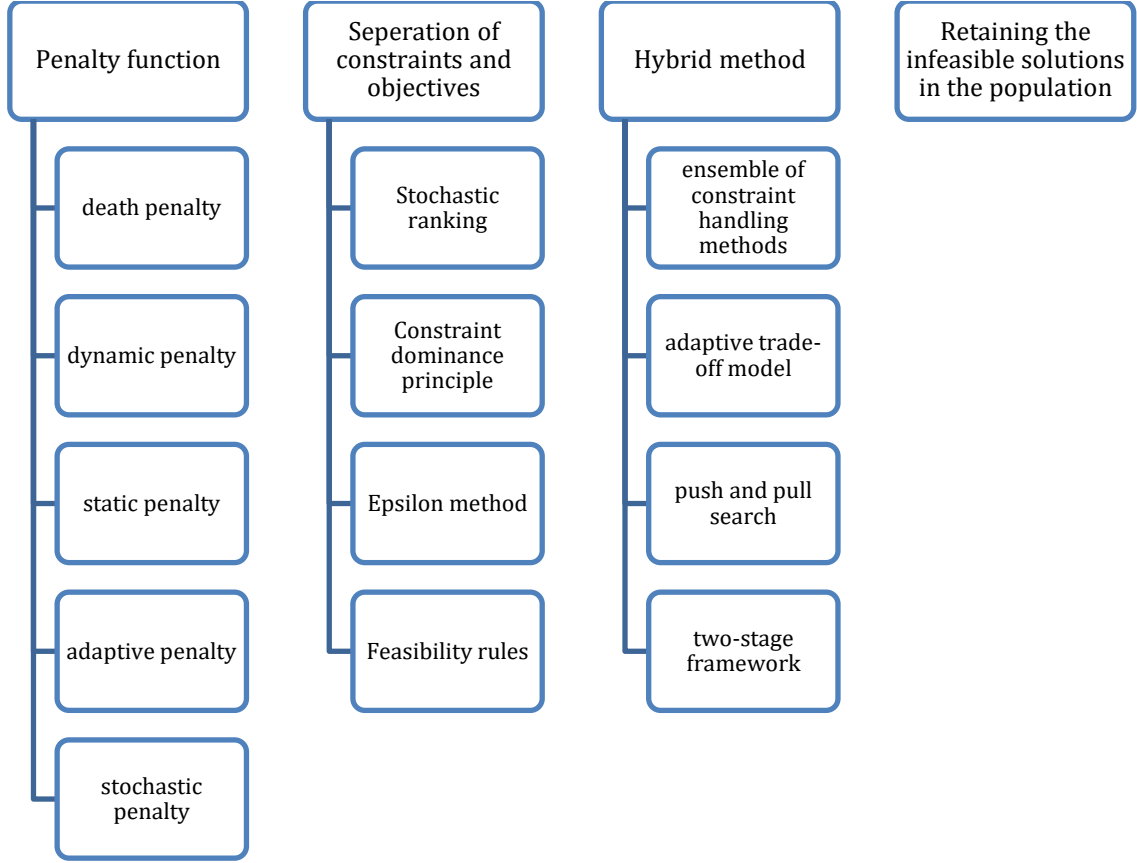

Figure 2 Taxonomy of different constraint handling methods in MOEA

Generally, penalty function techniques are one of the most simple CHTs. There are several types of penalty functions used with EAs, which the most important ones include [56]:
- Death penalty
- Dynamic penalty
- Static penalty
- Adaptive penalty
- Stochastic ranking

Details regarding the penalty function methods will be discussed in the next section.

## 3.1 Penalty function approach

One of the easiest and most common ways to handle constraints in multiobjective evolutionary algorithms is the penalty function method.

From a mathematical point-of-view, two types of penalty functions could be considered as follows:
- Interior methods
- Exterior methods

The first type of penalty functions, interior methods, the penalty factor is selected such that the value will be small away from the constraint boundaries and need an initial feasible solution [55], whereas exterior methods do not need an initial feasible solution [5]. Also, it should be noted that some of the infeasible solutions should be retained in the populations so that they are able to converge to a solution, which lies in the boundary between the feasible and infeasible regions [57].

In the penalty function method, any infeasible solution is ignored [58]. First, all constraints should be normalized, and for each solution, the constraint violations are calculated as follows (Equation 4):

$$w_j(x^i) = \begin{cases} |\bar{g}(x^i)|, & \text{if } \bar{g}(x^i) < 0 \\ 0, & \text{otherwise} \end{cases} \quad (4)$$

where $\bar{g}(x^i)$ refer to the normalized values for a given constraint $\bar{g}_j(x^i) \geq 0, j = 1, \dots, J$. Once the violations for the constraints are calculated, the values are added to determine the overall violation as follows (Equation 5):

$$\Omega(x^i) = j \quad (5)$$
$$= \sum_{j=1}^{J} w_j(x^i)$$

Also, a penalty parameter is multiplied to the sum of constraint violations and then added to the objective function values. If a proper penalty parameter is selected, MOEAs will work well; otherwise, a set of infeasible solutions or a poor distribution of solution is possible.

### 3.1.1 Static penalty functions

In static penalty proposed by [59], the penalty coefficient increases as a higher level of violation is reached. In fact, penalty functions do not change, a static penalty function is suggested, and several levels of violation are introduced in which the static penalty parameter is changed in case higher levels of violation are achieved [60]. In the static penalty function, the expanded objective function is (Equation 6):

$$\varphi(x) = f(x) + \sum_{j=1}^{p} C_{kj} G_j \quad (6)$$

Where $G_j = max\{0, g_j(x)\}^\beta$; and k=1,…,l where l presents the number of levels of violation.

### 3.1.2 Dynamic penalty functions

In this category, functions are changed based on the iteration number. [61] Proposed the following dynamic penalty function, in which the penalty increased when the iteration number increases.

Dynamic multiobjective optimization problems (DMOOPs) involve the simultaneous optimization of different objectives subject to a number of given constraints, where the objective functions, constraints, and/or dimensions of the objective space could change over time. EAs have acquired great attention among researchers for solving the above-mentioned problems.

### 3.1.3 Adaptive penalty functions

In this category, infeasible individuals are penalized according to the feedback taken from the search process [62]. [63] proposed a CHT based on the adaptive penalty function and distance measure, which both change as the objective function value and constraint violations of an individual varies.

Penalty-based constraint handing for multiobjective optimization is similar to single-objective problems in which a penalty factor is added to all the objectives. [63] Proposed a self-adaptive penalty function, which is suitable for solving constraint multiobjective optimization problems using evolutionary algorithms. In the self-adaptive penalty function method, the amount of penalty added to infeasible individuals are identified by tracking the number of feasible individuals. Also, the method uses improved objective values instead of the original objective function values [64].

### 3.1.4 Annealing-based penalty functions

[65] Introduced a multiplicative penalty function based on simulated annealing. In this type of penalty function, the temperature is decreased when the iteration number increases, which leads to an increased penalty.

### 3.1.5 Co-evolutionary-based penalty functions

[66] Proposed a co-evolutionary approach in which the population is partitioned into two subpopulations. The first population evolves solutions, and the second population evolves penalty factors. In this approach, the penalty function considers information taken from the amount of constraint violations and a number of violations.

There are other types of penalty function methods, which table 1 presents a summary and critique of the techniques for constraint handling.

Table 1 The important CHTs (penalty function)- *RQ2*

| Method | Criticism | Consequences |
|---|---|---|
| Death Penalty [12] | • No information is used from infeasible points.<br>• It may require initialization of the population and lack of diversity [67]. | • Consumes many evaluations<br>• Low success rate |
| Static Penalty [68] | • It is required to set up a high number of penalty parameters.<br>• It is also problem-dependent. | • Time consuming |
| Dynamic Penalty [69] | • It is hard to drive good dynamic penalty functions in real case.<br>• In some cases, this method converges to either an infeasible or feasible solution that is far from the global optimum [70]; [71] | • Premature convergence or even infeasible solution in some cases |
| Adaptive Penalty [72] | • Setting the parameters is difficult, such as determining the appropriate generational gap.<br>• It requires the definitions of additional parameters [73]. | • Time consuming |
| Annealing Penalties [74] | • The main disadvantage is its sensitivity to the values of its factors.<br>• To handle linear constraints, the user should provide an initial feasible point to the algorithm. | • The performance of the algorithms is not good. |
| Self-adaptive Penalty [75][76] | • It defines four additional parameters that may affect the fitness function evaluations. | • Time consuming & weak or strong penalty during evolution |
| Segregated genetic algorithm (SGA) [77] | • The main difficulty is selecting the penalties for each of the two sub-populations. | • Time consuming |
| Penalty function based on feasibility [13] | • The main issue is maintaining diversity in the population, and in some cases, use of a niching method combined with higher-than-usual mutation rates is essential. | • Premature convergence |

### 3.2 Separation of objective function and constraints

Unlike the penalty function technique, another approach exists that separates the values of objective functions and constraints in the nature-inspired algorithms (NIAs) [78], which is known as separation of objective function and constraints. [79] Initially proposed the idea to divide the search space into two phases. In the first phase, feasible solutions, are found, and in the second phase, optimizing the objective function is considered.

Representative methods of this type of CHT are as follows:
- Stochastic ranking (SR)
- Constraint dominance principle (CDP)

- Epsilon CHT
- Feasibility rules

The next section provides further details about this type of constraint handling method.

### 3.2.1 Constraint dominance principle

In constraint dominance principle (CDP), three feasibility rules are applied to compare any two solutions. If $x^1$ is feasible and $x^2$ is infeasible, then $x^1$ would be better than $x^2$. If both solutions are infeasible, then the solution with a smaller constraint violation is better. If both are feasible, then the one dominating the other is better. [16] Adopted CDP to handle constraints in NSGAII (NSGAII-CDP), in which the population is divided into feasible and infeasible sub-populations. NSGAII-CDP first selects offspring from the feasible solutions and then selects solutions from the infeasible solutions.[80] Also adopted CDP to handle constraints in the MOEA/D framework.

### 3.2.2 ε-constrained (EC) Method

The basic principle of the ε-constrained method, first introduced by [81], is similar to the superiority of feasible solution (SF) proposed by [82] (Eqs. 12-13). The epsilon value is updated until the parameter k reaches the control generations $T_c$. [83] Embedded the epsilon CHT in MOEA/D so that the epsilon value is set adaptively of r comparison. Also, the violation threshold is based on the type of constraints, size of the feasible space, and the search outcome. In the method proposed by [83], the infeasible solutions with violations less than threshold are identified (Eqs. 7-8).

$$\varepsilon(0) = V(x_\theta) \quad (7)$$

$$\varepsilon(k) = \begin{cases} \varepsilon(0)(1 - \frac{k}{T_c})^{cp}, & 0 < k < T_c \\ 0, & k > T_c \end{cases} \quad (8)$$

where $x_\theta$ presents the top $\theta$th individual at initialization; and the cp parameter is selected between [2,10] [84].

### 3.2.3 Feasibility rules

The popularity of this method depends on its ability to be coupled to a range of algorithms without announcing new parameters (factors) [30].

The feasibility rules proposed by [13] are simple, could be integrated into a variety of algorithms without adding new parameters, and thus, are largely used in the research field. [85]–[87] developed feasibility rules for the selection process, which have been adopted by different evolutionary algorithms such as DE, PSO, and GA. According to the number of feasible solutions, the search space could be divided into three phases as follows [84]:

- No feasible solution is found.
- There exists at least one feasible solution.
- Integrating the parent-offspring population has more feasible solutions than the size of the next generation population.

The feasibility rules used in multiobjective optimization, also known as the superiority of feasible solution (SF), are addressed as follows [84] (Equation 9):

$$fitness_m(x) = \begin{cases} f_m(x) & \text{if } x \text{ is feasible} \\ f_{worst}^m + v(x) & \text{otherwise} \end{cases} \quad (9)$$

where $f_{worst}^m$ and v(x) show the mth objective value of the worst feasible solution and the overall constraint violation, respectively.

#### 3.2.3.1 Feasibility rules in Differential Evolution (DE)

Although the feasibility rules introduced by [13] have also been widely used by other researchers in DE [85], [86], [88]–[94], they have been rarely used in multiobjective differential evolution. Particularly, [95] used Pareto dominance in constrained space instead of the sum of constraint violations. Later, [96] adopted the Pareto dominance in Generalized Differential Evolution (GDE), but encountered difficulties when there exist more than three constraints and/or objective functions.

[97] Proposed a scheme for partitioning the objective space using the conflict information for multiobjective optimization. [98] introduced an operational efficient model based on Data Envelopment Analysis (DEA) and introduced DE along with the feasibility rules to optimize the mentioned model. [99] proposed a combined constraint handling framework, known as CCHF, for solving constrained optimization problems, in which the features of two well-known CHTs (i.e. feasibility rules and multiobjective optimization) were addressed in three different situations (feasible situation, infeasible situation, and semi-feasible situation).

.

#### 3.2.3.2 Feasibility rules in PSO

[100] Employed feasibility rules as a constraint handling technique to recognize the most competitive PSO variant when solving constrained numerical optimization problems (CNOPs). In the research by [100], local-best was identified to be better than global best PSO. [101] [102] adapted an artificial bee colony algorithm (ABC) to solve CNOPs by using feasibility rules by modifying the probability assignment for the roulette wheel selection.[103] Compared different GA variants using the feasibility rules as the constraint handling method. [104] Introduced a hybrid version of PSO to solve constrained optimization problems and found that the swarm at each generation is split into several sub-swarms. Also, the hybrid version applied the feasibility rules to compare particles in the swarm.

#### 3.2.3.3 Feasibility rules in GA

[105] Proposed a GA with a new multi-parent crossover for solving constraint optimization problems, in which the feasibility rules were added to handle the constraints. The latter authors also solved constraint numerical optimization problems by using different GA variants along with feasibility rules and found that all GAs perform equally. [22] Introduced a two-phase framework for solving constraint optimization problems. Specifically, the first phase ignores the objective function and the genetic algorithm minimizes the violation of the solutions, while the second phase optimizes bi-objective functions, including the original objective function and constraints satisfaction. Moreover, feasibility rules is applied to assign fitness values to the individuals.

#### 3.2.3.4 Feasibility rules in other population-based algorithms

Feasibility rules have been adapted to other population-based algorithms, such as artificial immune systems [106]–[111], organizational evolutionary algorithm [112], biogeography-based optimization [113], and bacterial foraging optimization [114].

### 3.3 Retaining infeasible solutions in the population

Another CHT is used to retain the infeasible individuals in the population. In other words, a constraint multiobjective optimization problem with m objective is transformed to an optimization problem with m+1 objectives, which could save the infeasible solution during the evolution process [28], [115]. [116] proposed a constraint handling technique so that the individuals that possess low Pareto rank and low constraint violation will be chosen.

### 3.4 Hybrid methods

Hybrid methods combine several CHTs to handle constraints. [117] addressed four different hybrid methods in their work: 1) the ensemble of constraint handling method [84] separates the population into three sub-populations; 2) the adaptive trade-off model [118] contains two different CHTs; 3) in the push and pull search [119], the population is pushed to the unconstrained Pareto front (push) then the population is pulled back to the Pareto front (pull); 4) the two-phase framework (ToP) [120] solves a constraint multiobjective optimization problem by first converting the objective functions into a single objective function via the weighting method then, in the second phase, a constrained MOEA is adopted to attain the Pareto feasible solutions. The disadvantages for each category in Figure 1 are summarized in Figure3.

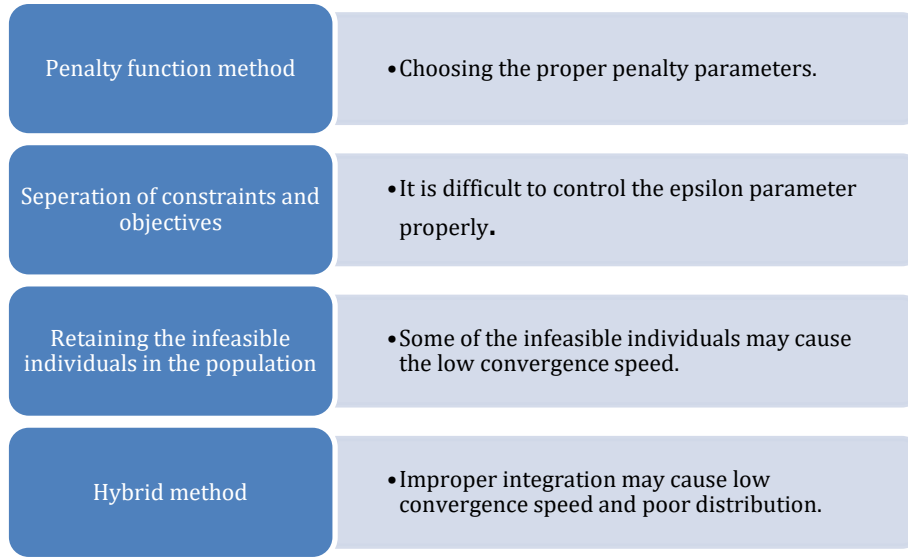

Figure 3 The main disadvantages of CHTs in MOEA- RQ2

The authors of [30] Presented a state-of-the-art taxonomy of CHTs, which is illustrated in Figure 4. Figure 5 presents the different state-of-the-art CHTs that have been used since 2000. As mentioned, [13] first applied feasibility rules to the genetic algorithm. [121] introduced stochastic ranking, which employs a user-defined parameter instead of using penalty factors and is able to control the infeasible solutions based on the sum of constraint violation and objective function values. [33] Proposed the epsilon-constraint method that transforms the constraint optimization problem into an unconstrained problem. [33] Addressed a multi-constrained optimization problem based on the KS function. [122] proposed a boundary search approach inspired by the ant colony metaphor based on conducting a boundary search between a feasible and infeasible solution. [28] Proposed an additional objective to solve a bi-objective optimization problem, where the first objective is the original problem and the second objective is the constraint violation measure. [36] Addressed a combination of four CHTs, namely feasibility rules, stochastic ranking, self-adaptive penalty function and the epsilon-constraint method to solve constraint numerical optimization problems. The first and second techniques are discussed in detail in the next sections, while the other methods presented in Figure 3 have been addressed previously in this study.

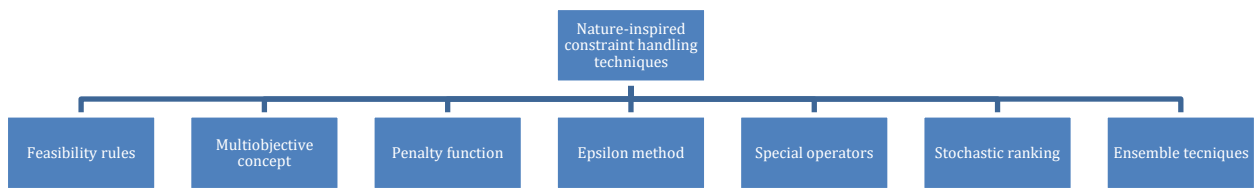

Figure 4 State-of-the-art CHTs [30]

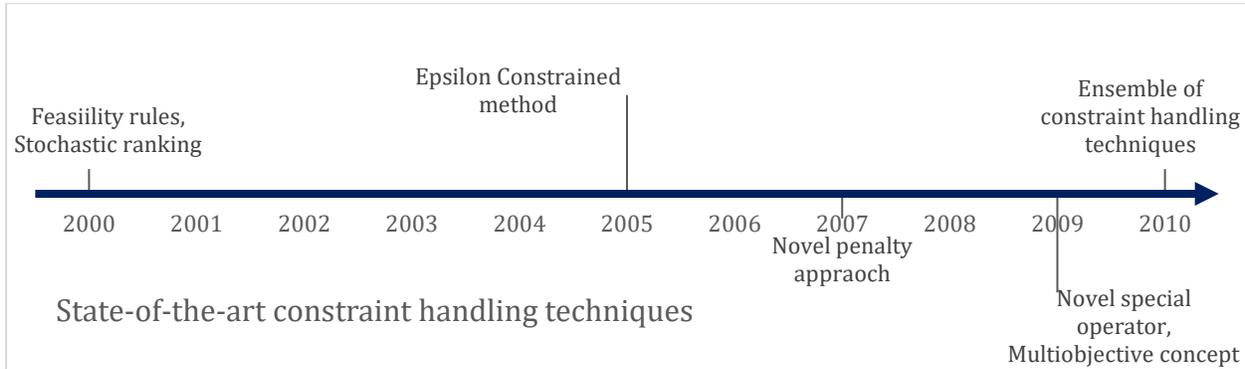

Figure 5 Timeline of different state-of-the-art CHTs

### 3.5 Stochastic ranking

[121] proposed the stochastic ranking (SR) approach to balance between the objective and penalty functions stochastically. The method was tested using a strategy evolution on several benchmarks, and the results showed that the method is able to improve the search performance with a user-defined parameter without introducing complicated variation operators. SR has also been coupled with other population-based algorithms, such as ant colony optimization (ACO) [123], [124], differential evolutionary (DE) [125]–[127], and evolutionary programming (EP) [128].

### 3.6 Ensemble techniques

Ensemble CHTs provide a new research platform to tackle constrained multiobjective optimization problems. Combining several CHTs could improve the capability of an approach compared with a single CHTs [30], [129]. For instance, [130] proposed a combination of four CHTs, namely nondominated sorting, constrained-domination principle, multiple constraint ranking, and dynamic penalty function, and incorporated the proposed technique into an MOEA based on NSGAII. Some other ensemble CHTs have been reported [36], [131], [132]. Although the ensemble CHT has a competetive performance, it suffers from being parameter-dependent.

### 3.7 Multiobjective concept

Based on the multiobjective optimization concept, a constraint single-objective optimization problem is transferred to an unconstrained multiobjective optimization problem [5]. The multiobjective version of the optimization problem possesses an extra objective function, which presents the sum of constraint violation [28], [133]–[135].

The authors of [136] presented a taxonomy for constraint handling strategies in multiobjective GA, which include:
- Penalty function methods
- Separation method
- Special operators
- Repair methods

Among these strategies, the penalty function method is not straightforward in multiobjective GA since the fitness assignment is based on the non-dominance rank of a solution rather than its objective function values [137]. Yet, the penalty function method is one of the most popular CHTs in constraint multiobjective optimization. Whenever the multiobjective function and constraint violation for each constraint are assessed, the sum of violations is added to each objective function value considering the multiplication of the penalty parameter [138], [139].[19] Proposed two approaches, namely OEGADO and OSGADO. OEGADO runs several GAs in parallel so that each GA optimizes one objective, whereas OSGADO runs each objective sequentially with a common population for all objectives.

### 3.8 Repair approaches

There are several techniques used as repair algorithms, in which the search space is reduced (since only feasible individuals are considered):
- In the permutation encodings method, each solution of an EA population is simply signified as an ordered list [140], [141].
- Repair procedures in binary representations, which could be shown as fixing the number of 1s in binary representations and Hopfield networks [142], [143].
- Repair methods in graphs that are represented as spanning trees and repairing graphs [144], [145].
- Repair methods in grouping GAs, which are proper for scenarios that a number of items should be assigned to a set of groups [146], [147].

Pure EAs do not perform well in complex combinatorial problems with a high number of constraints [148], [149]. Single-solution based algorithms (e.g. local search, simulated annealing) have good performance in exploitation, while population-based algorithms (e.g. swarm intelligence, EA are exploration-oriented. In these problems, hybridization of population-based algorithms with single-based algorithms can improve the power of both exploration and exploitation [150]–[152]. A memetic algorithm is a hybridization of an EA and a local search (LS) approach that LS is applied to improve the quality of the fitness function. On the other hand, LS could be used as a CHT [152], i.e. the local repair algorithm only consider feasible individuals leading to reducing the search space. Repair methods could be applied to EAs in several ways, such as in permutation encodings [153], [154], in binary representation [155], and in graphs and trees [156], [157]. Although repair algorithms have numerous advantages, some disadvantages do exist. For instance, repair algorithms are problem-specific and must be designed for a specific problem [158]. Table 2 shows a summary of the disadvantages of the state-of-art CHTs.

Table 2 A summary of disadvantages of the state-of-art CHTs- RQ2

| Method | Disadvantages |
|---|---|
| Ensemble method | Although the ensemble CHT has a competitive performance, the method is parameter-dependent. |
| Repair method | Repair algorithms are problem-specific and, thus, must be designed for a specific problem. |
| Feasibility rules | The method is likely to lead to premature convergence. |
| Stochastic ranking | Although the method has been employed in several nature-inspired algorithms, it is not often used for the multiobjective version of the algorithms. |
| Epsilon-constraint method | In some cases, premature convergence has been reported, while other works report that the method relies on gradient-based mutation. |
| Multiobjective concept | It may require gradient calculation [30]. |

## 4 Other approaches

Table 3 provide a summary of novel approaches proposed between 2020 and 2021 to tackle constrained multiobjective optimization problems. Based on Table 3, there is signs of a renewed interest in constrained multi-objective optimization, even the clear superior amount of research in constrained single-objective optimization.

Table 3 Novel approaches for constrained multiobjective optimization problems between 2020 and 2021

| Method | Source | Method | Source |
|---|---|---|---|
| KKT points for constrained multiobjective optimization | [160], [161] | Surrogate-assisted evolutionary algorithm | [159] |
| IoT and cloud computing | [163] | Purpose-directed two-phase multiobjective differential evolution | [162] |
| Indicator-based constrained handling technique | [165] | Directed Weight Vectors | [164] |
| Decomposition-based algorithm | [117], [167] | Gradient-based repair method | [166] |
| Push and pull search embeded | [169] | Detect and scape strategy | [168] |
| Multi-stage evolutionary algorithm | [171], [172] | Reference points-based method | [170] |
| Partition selection | [174] | multi-objective wireless network optimization using the genetic algorithm | [173] |

As a general, a taxonomy of CHTs in MOEAs could be summaraized in Figure 6. The CHTs presented in Figure 6 have been explained in detailes in previous sections. As it is mentioned earlier, most of constraint handling techniques developed for single-objective optimization problems can be applied to multi-objective optimization problems. It is worthy to note that among them, stochastic ranking [175][176] [177] [178] [179] [180], penalty function [181][63], multi-objective method [182] [183] [184] [185], Epsilon constrained method [186] [187] [188] [189] [190] [191] [192] [193], transforming method [194] [195] [196] [197] [198] [199] [200] [201] [202], feasibility rules (with a modification) [203][204], hybrid methods [205] [206] [207] [208] [209], and repair operators [210] [211] [212] [213] [214] [215] have been addressed to multi-objective optimization problems.

## 5  Benchmark test problems

To measure the performance of the evolutionary algorithms, many benchmark or test problems have been suggested. On the other hand, the benchmark problems help researchers to better understand the strengths and weaknesses of an algorithm. These test problems are classified as single-objective such as Rosenbrock [216], G01-G09 [217]
,Himmelblau's problem [218], Welded Beam [219], Pressure Vessel [220], Tension–Compression Spring [221],
Speed Reducer [222], Corrugated bulkheads design [223], Heater exchanger [224], Multiple disk clutch brake [225],
Rolling element bearing [226], Car side design [227], Stepped beam design problem [228], mult-objective including BNH [229], OSY [230], , ZDT [231], BT [232], Truss2D [233], and many-objective optimization problems for example C-DTLZ [234], WFG[235] , DTLZ [236]. Among the above-mentioned test problems, some them are still unconstrained. For more details, it is suggested to the review paper in the field by [235].

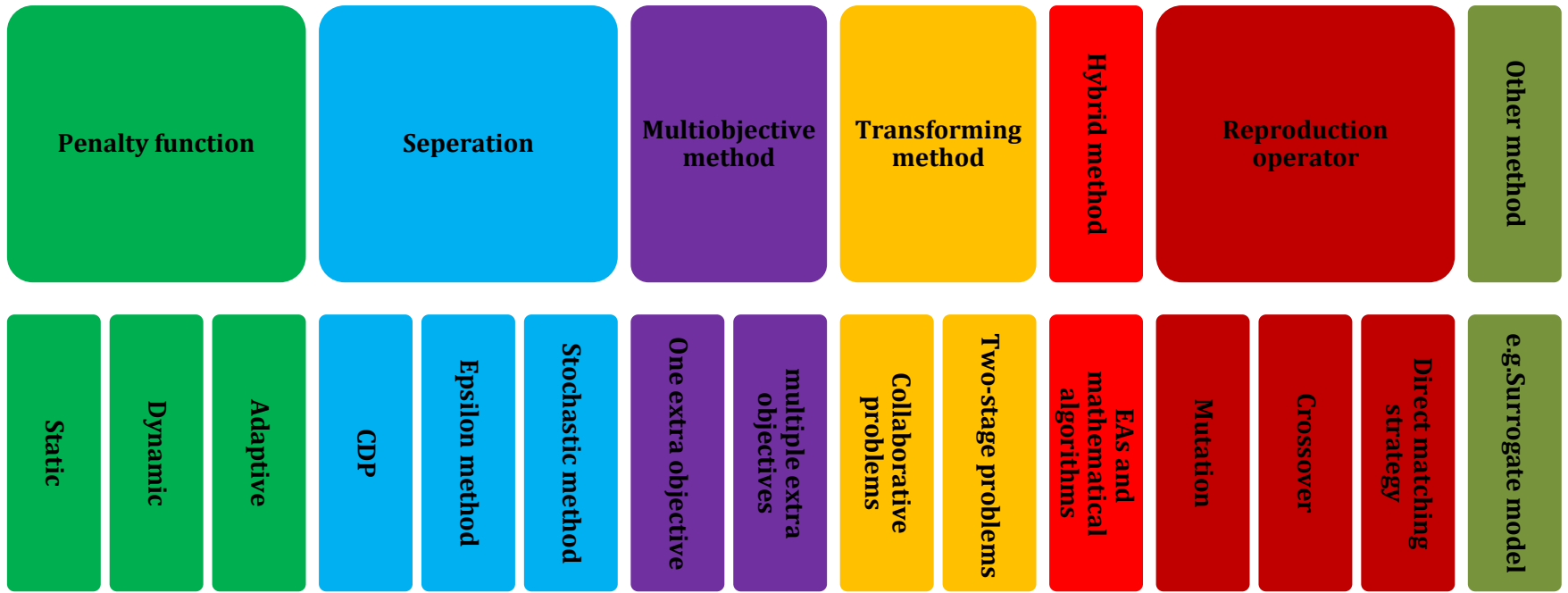

Figure 6 A generalized taxonomy of CHTS for multiobjective optimization problems

## 6 Scientometric Analysis (RQ3- RQ6)

Scientometric analysis is conducted to scientifically measure and analyze the literature in a particular field of study [237]. Bibliometrics is the most famous field of scientometrics that uses statistics to analyze and measure the impacts of books, research articles, conference papers, etc. [238]. Recently, this field of analysis has attracted much attention from researchers and has been used in various literature review fields [239]–[243]. To perform the analysis in this work, VOSviewer [46] and RStudio were used. The following sub-sections provide new insight into the scientometric analysis in the field.

### 6.1 Citation statistics

Figure 7 displays the trend of published documents, which shows that the number of documents in the field significantly increased from 2003 until the end 2020 (just above 100 documents).

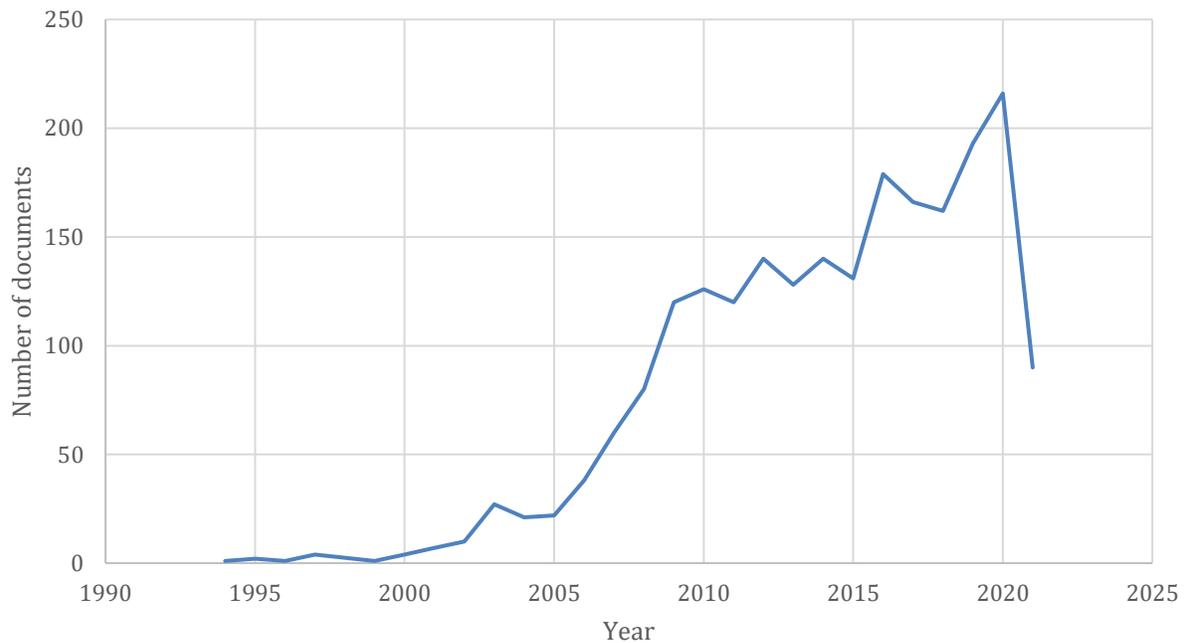

Figure 7 Trend of published documents

Figure 8 presents combo chart of the number of documents vs. total citations since 2008. In 2020, the most documents and citations were achieved (216 and 7544, respectively). It is apparent that the number of citations has increased dramatically according to the trend. According to the WOS, the numer of citations of the top articles in the field was analyzed and is presented in Supplemenatary C (Figure 1). Of the 735 related documents in WOS, about 45824 citations were identifed from the related papers with an average of 1992.35 citations per year and an average of 62.35 citations per item. [16], [244], and [245] are top 3 cited articles with 20013, 2609, and 1591 citations in WOS, respectively.

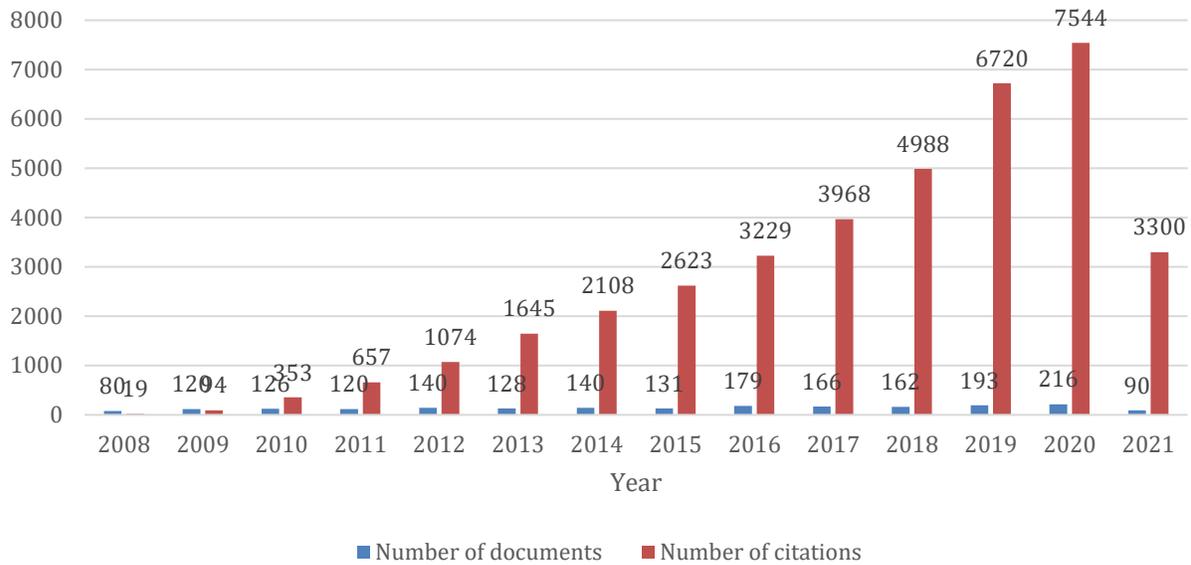

Figure 8 Combo chart of number of documents vs. total citations (Scopus)

### 6.2 Statistics based on document types

Among the document types, including article, proceedings paper, review, and other items indexed by WOS, a total of 735 publications on constraint handling multiobjective population-based optimization algorithms was found (Table 4). From the search, articles were the most popular document type, comprising a total of 522 articles (71.02% of 735 documents) with 2.60 authors per publication (APP). Also, articles as the document type had the highest CPP2020 of 84.10, followed by proceedings papers with TP of 220 (29.93% of contributions and APP=2.13). Moreover, there is a significant difference between the TC2020 of article and that of proceedings paper.

Figure 9 presents the distribution of documents based on different types, according to WOS. It is clear from the figure that conference papers have the most contributions before 2010 followed by articles. However, since 2010, articles have the most contributions in the field. It is also interesting to note that book/book chapters have been published since 2000, however, the most number of book/book chapters have been published after 2010.

Table 4 Citations analysis based on document type

| Document type | TP | % | AU | APP | TC2020 | CPP2020 |
|---|---|---|---|---|---|---|
| Article | 522 | 71.02 | 1362 | 2.60 | 43,904 | 84.10 |
| Proceedings paper | 220 | 29.93 | 469 | 2.13 | 1,543 | 7.01 |
| Review | 16 | 2.17 | 20 | 1.25 | 806 | 50.37 |
| Other items | 23 | 3.12 | 134 | 5.82 | 468 | 20.34 |

TP, AU, APP, TC2020, and CPP2020 present total number of articles; total number of authors; total number of authors for each publication; total citations from WOS since publication year to the end of 2020; total citations for each paper, respectively; Other items: early access and letters [246].

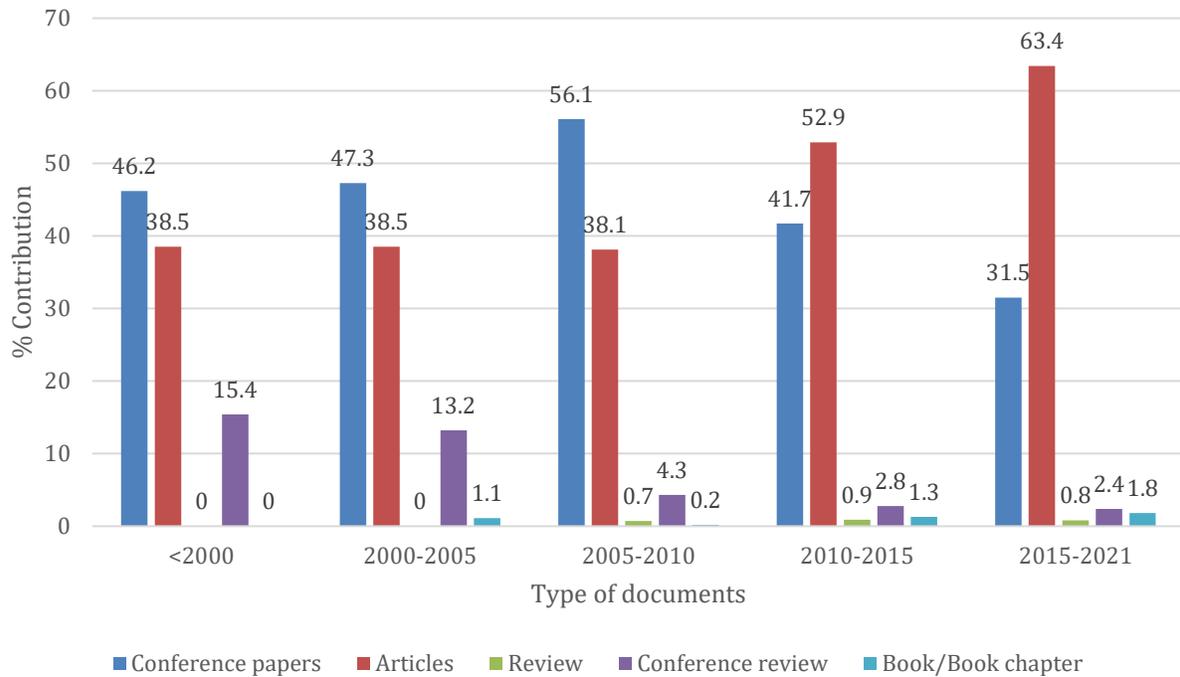

Figure 9 Type of research outputs

### 6.3 Publication statistics based on journal

Table 5 presents the top 20 journals that have published the greatest number of constraint handling multiobjective population-based algorithms papers based on Scopus. Accordingly, Lecture Notes In Computer Science (117), Applied Soft Computing Journal (57), and Swarm and Evolutionary Computation (30) are most are the most utilized, which predominate in the field of optimization and evolutionary computations.

A total of 735 articles were published in 399 journals, which are classified among the 51 WOS categories in SCI-EXPANDED. Table 6 lists the 10 most productive WOS categories. A total of 271 articles (36.87% of 735 articles) were published in the first category (Computer Science Artificial Intelligence), of which 83.39% were published in Engineering Electrical Electronic (23.26%) and Computer Science Theory Methods (23.26%). Comparing the top 10 categories, the highest CPP2020 of articles published in the Computer Science Theory Methods category is 190.011, which includes the paper entitled: "A fast and elitist multiobjective genetic algorithm: NSGA-II" by [16], and the highest APP for articles published in the Energy Fuels category is 2.97.

Table 5 The top 10 sources that have published the greatest number of constraint handling on GA papers (Scopus)

| # | Scopus | # of Documents | # | Scopus | # of documents |
|---|---|---|---|---|---|
| 1 | Lecture Notes In Computer | 117 | 11 | IEEE Access | 32 |
| 2 | Applied Soft Computing Journal | 57 | 12 | Swarm and Evolutionary | 30 |
| 3 | "International Journal of Electrical | 27 | 13 | Engineering Optimization | 21 |

| 4 | "Kongzhi Yu Juece Control And | 13 | 14 | Soft Computing | 16 |
| 5 | Energy Conversion And | 12 | 15 | Studies In Computational | 16 |
| 6 | IEEE Transactions On Cybernetics | 12 | 16 | "Advances In Intelligent Systems | 15 |
| 7 | Structural And Multidisciplinary | 13 | 17 | "Communications In Computer | 14 |
| 8 | IEEE Transactions On | 27 | 18 | Engineering Applications Of | 14 |
| 9 | Electric Power Systems Research | 10 | 19 | Energy | 13 |
| 10 | Applied Intelligence | 12 | 20 | Information Sciences | 13 |

Table 6 The top 10 productive WOS categories

| # | Web of Science category | TP | AU | APP | TC2020 | CPP2020 |
|---|---|---|---|---|---|---|
| 1 | "Computer Science Artificial Intelligence" | 271 | 628 | 2.31 | 35,074 | 129.42 |
| 2 | "Engineering Electrical Electronic" | 171 | 463 | 2.70 | 3,688 | 21.56 |
| 3 | "Computer Science Interdisciplinary Applications" | 92 | 243 | 2.64 | 2,691 | 29.25 |
| 4 | "Operations Research Management Science" | 61 | 131 | 2.14 | 1,541 | 25.26 |
| 5 | "Computer Science Theory Methods" | 171 | 385 | 2.25 | 32,492 | 190.011 |
| 6 | "Engineering Multidisciplinary" | 64 | 167 | 2.60 | 2,377 | 37.14 |
| 7 | "Mathematical interdisciplinary applications" | 45 | 116 | 2.57 | 645 | 14.33 |
| 8 | "Energy fuels" | 41 | 122 | 2.97 | 950 | 23.17 |
| 9 | "Computer Science Information Systems" | 48 | 124 | 2.58 | 746 | 15.54 |
| 10 | "Automation Control Systems" | 60 | 155 | 2.58 | 1,178 | 19.63 |

TP, AU, APP, TC2020, and CPP2020 present total number of articles; total number of authors; total number of authors for each publication; total citations from WOS since publication year to the end of 2020; total citations for each paper, respectively ; Other items: early access and letters [247]

Figure 10 provides a comparison of the development trends of the top four productive WOS categories, including "Computer Science Artificial Intelligence", "Engineering Electrical Electronic", "Computer Science Theory Methods", and "Computer Science Interdisciplinary Applications".

Between 2001 and 2021, Computer Science Artificial Intelligence was the most predominant category and has possessed the highest number of publications since 2004, excluding the period between 2007 and 2008. The three other categories possess fluctuations between 2001 and 2021, and as of writing this paper, "Computer Science Theory Methods" and "Engineering Electrical Electronic" have the same TP of 171.

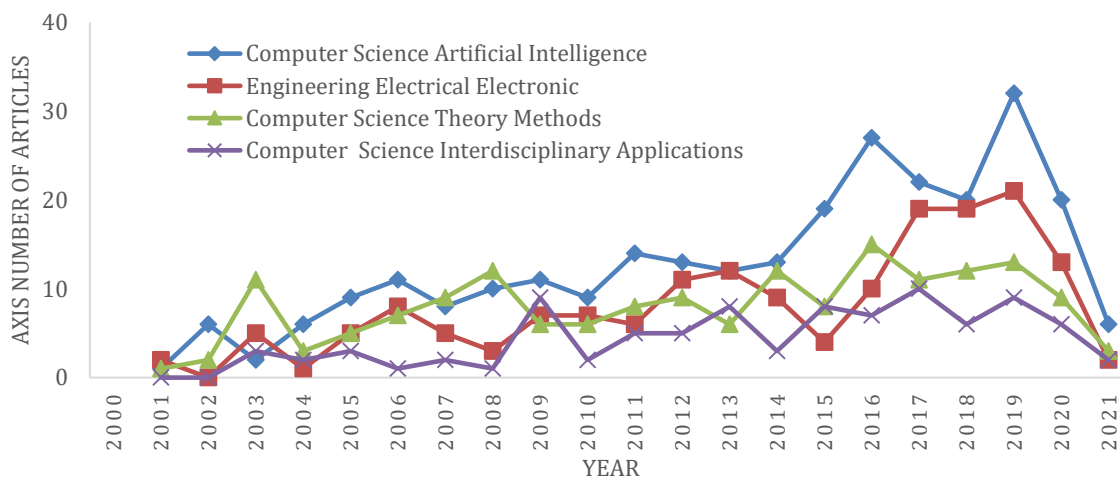

Figure 10 Comparison of the development trends of the top four productive WOS categories

### 6.4 *Publication statistics by countries*

Figure 11 presents the distribution of documents by the top most active countries in both databases. It is apparent that China, India, and the USA are the top three active countries in the field according to Scopus (respectively), while China, the USA, and India are the top 3 active territories in the field based on WOS, respectively. It is pertinent to mention that the USA is ranked second based on WOS, but India is ranked second according to Scopus. Also, it can be seen that there is a significant difference between the first rank (China) and second rank (India) based on the number of publications indexed by Scopus. Moreover, Figure 12 presents the collaboration among countries, where the links across the circles depict the collaborations, and the size of the circles represents the activities of the countries in the field. The green and yellow colors present the keywords that have been used recently, while the dark blue color indicates those that were used earlier (around 2008).

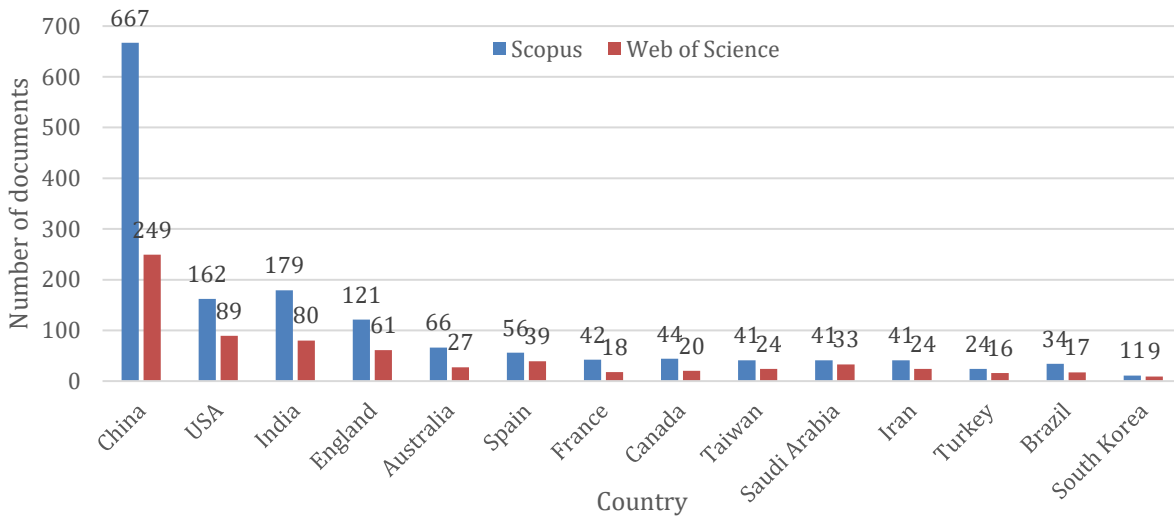

Figure 11  Research output of top 10 most productive countries across all databases

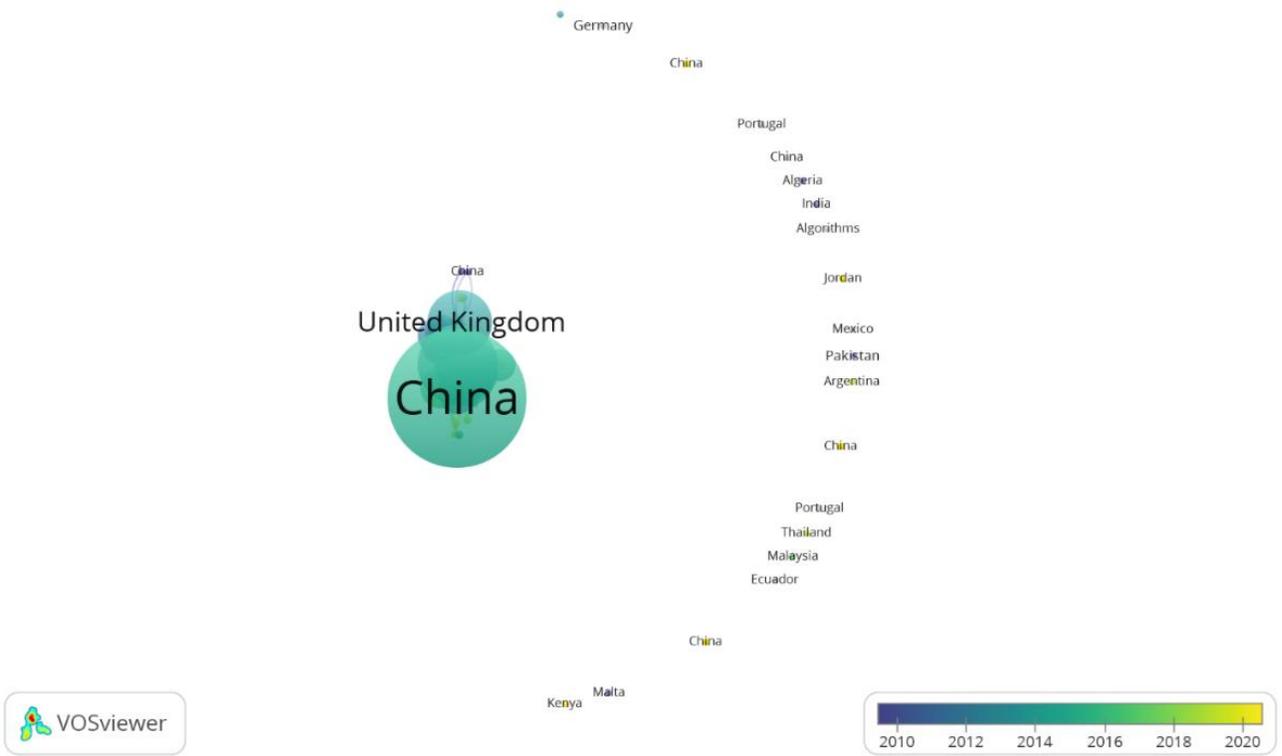

Figure 12 Country network visualization

Figure 13 displays the growth rate of the top 5 active countries in comparison to the world. While China and the USA have smooth trends between 2000 and 2020, India, the UK, and Australia show some fluctuations. Between 2002 and 2003, India presents the highest growth rate, then the trend continues smoothly until 2014 when it increases until 2015. The trend for the UK shows two growths between 2002-2003 and 2007-2008. While the number of articles published by Australia is much less than the four other countries, there is a significant rise between 2015 and 2016.

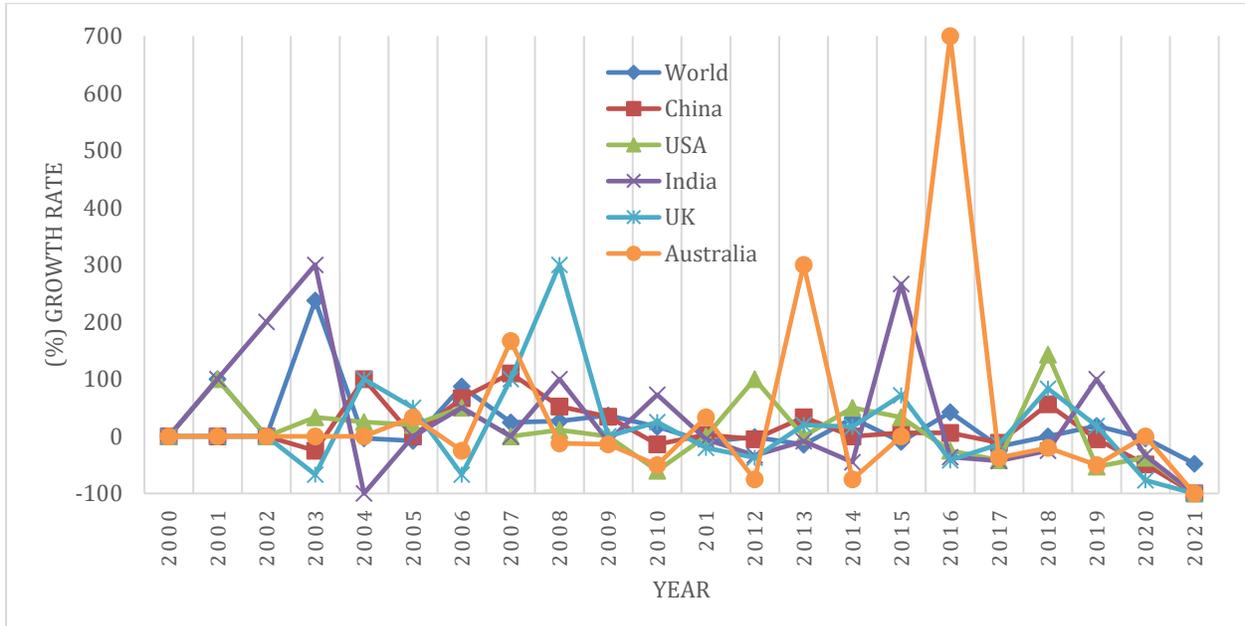

Figure 13 Growth rate of published documents for top 5 countries

### 6.5 Statistics based on the subject area

Figure 14 presents the distribution of articles based on the subject area. Computer science, Engineering, and Mathematics possess the most contributions with 936, 619, and 580 published articles, respectively. Comparatively, Pharmacology, Medicine, and Economics own the least contributions with 1, 4, and 6 published documents in the field, respectively.

Figure 14 Distribution based on the subject area (Scopus)

## 6.6 Statistics based on authors

Figures 15 shows the top authors with the most publication according to Scopus, (Supplementary C, figure 2 presents statistics based on WOS). Kalyanmoy Deb from "Michigan State University (USA)", Ray T. from "University of New South Wales (Australia)", and Carlos A. Coello Coello from "Cinvestav-IPN (Mexico)" with 38, 32, and 28 publications are the top 3 authors in the field as indexed by Scopus. WANG Y from "City University Hong Kong (Hong Kong)", Carlos A. Coello Coello from "Cinvestav-IPN (Mexico)", and Ray T. from "University of New South Wales (Australia)" are the top 3 authors in the area with 21, 20, and 17 documents (indexed by WOS), respectively. According to WOS, 1717 authors have worked on constraint multiobjective population-based optimization.

In total, 0.5241% of authors own more than 10 documents; 1.6307% possess between 5 and 10 documents; 4.5428% have between 3 and 5 papers; 11.7647% own 2 papers; and 81.5377% possess 1 document (Figure 16). Figure 17 presents the collaboration among the authors, where links across the circles depict the collaborations, and the size of the circles shows the activities of the authors in the field. In addition, the yellow color represents recent activity, and the dark blue color depicts the contributions prior to 2014.

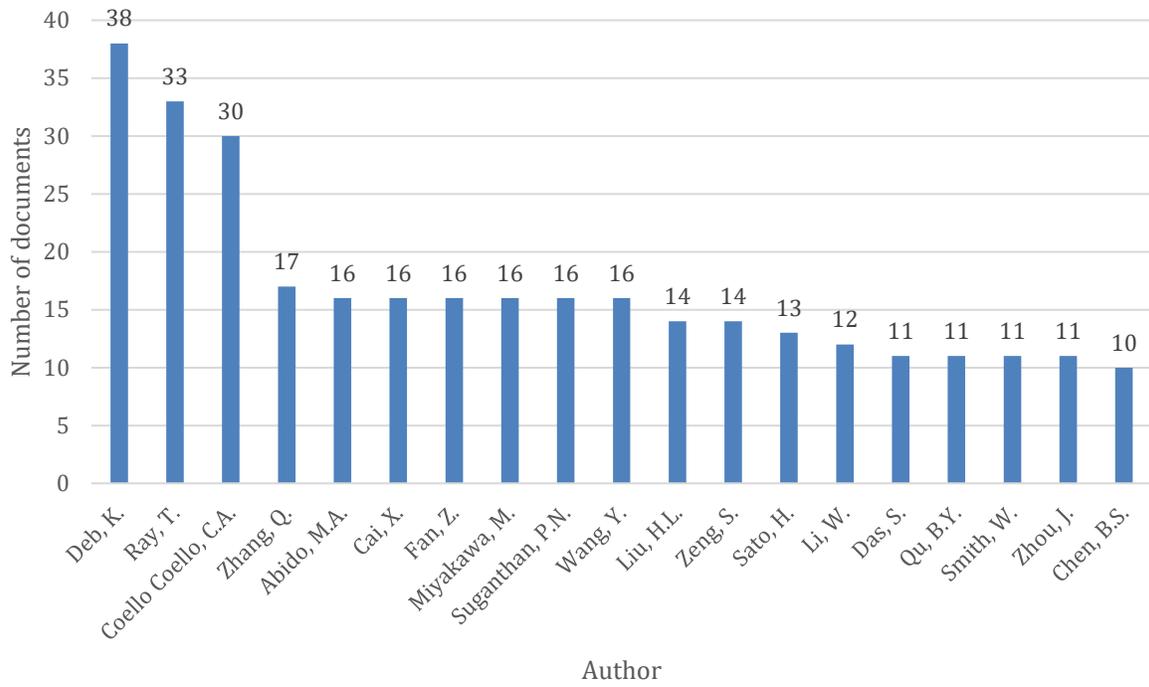

Figure 15 The most active authors in the field (Scopus)

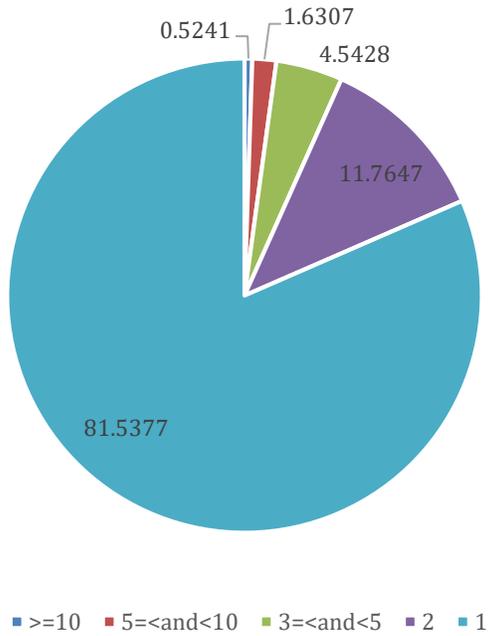

Figure 16 Contribution of authors based on the number of documents

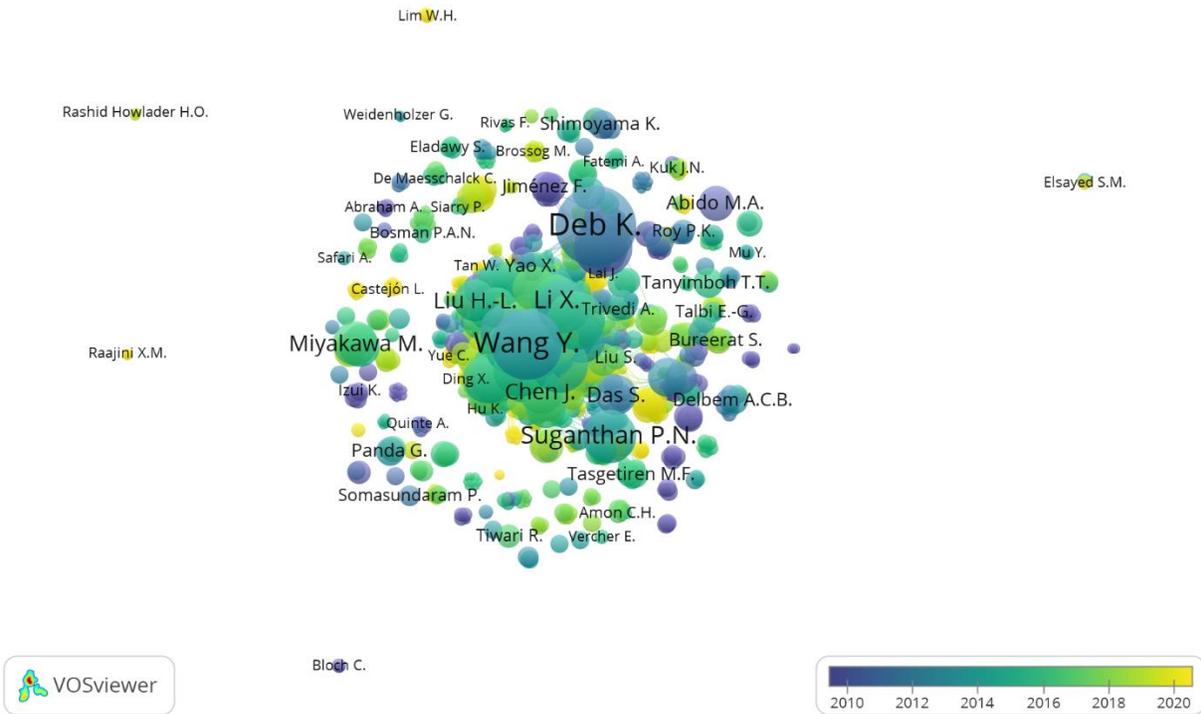

Figure 17 Collaboration among the authors (overlay visualization)

### *6.7 Statistics on keywords*

Keywords indicate the basic parts of a certain field of research and could offer insight into the organization and knowledge provided in the articles. Figure 18 provides an overlay visualization of the co-occurrence analyses via a network map based on the Scopus database. Each node in the network represents a keyword, and the link between nodes indicates the co-occurrence of the keywords. The top keywords in Scopus include multiobjective optimization, Pareto optimal solution, evolutionary algorithm, artificial intelligence, machine design, stochastic systems, distributed power generation, and reliability. The color of each circle represents the identified cluster, and the size of each circle illustrates the importance of the keywords; in other words, keywords with a larger circle have been used more than others. The green and yellow colors show the keywords that have been used recently, while dark blue color are those that have were used earlier (around 2008). Tables 7 presents the top keywords of 1-word, 2-word, and 3-word lengths extracted from Scopus. Specifically, optimization, algorithm, and scheduling are the top 1-word length keywords indexed by Scopus; genetic algorithm, constraint handling, and constrained optimization are the top 2-word length keywords; and constraint-handling techniques, Particle Swarm Optimization (PSO), and multiobjective optimization are the top 3-word length keywords indexed by Scopus.

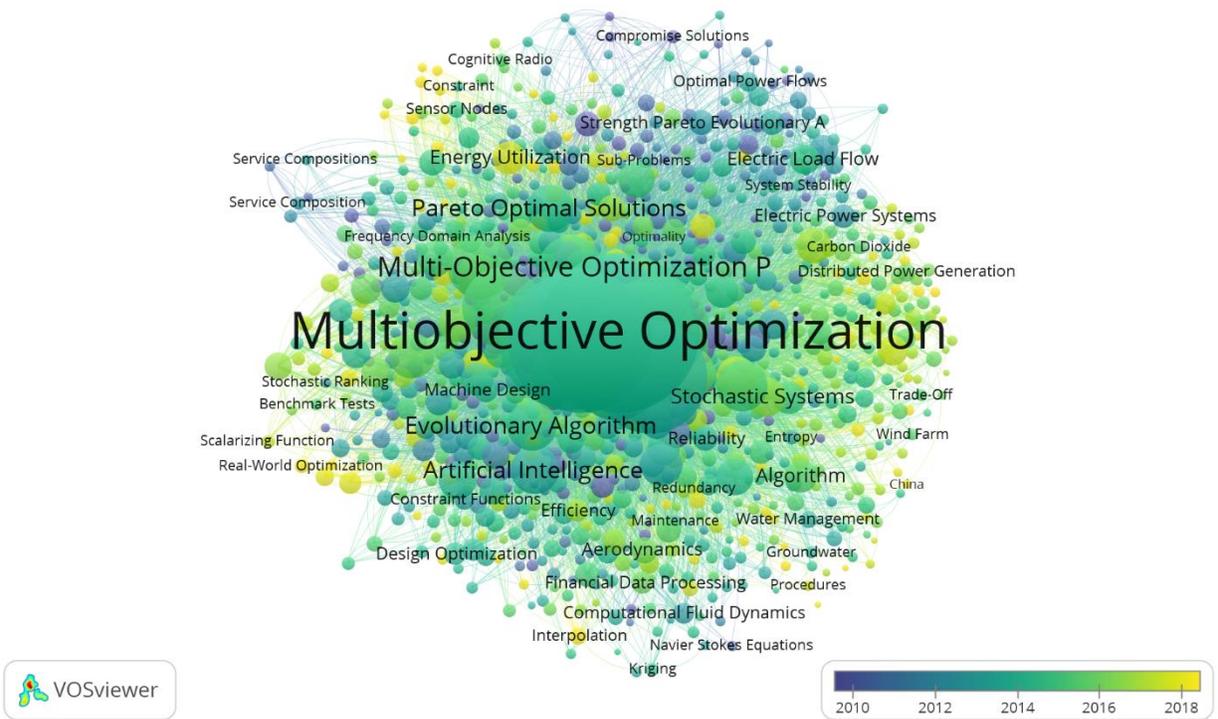

Figure 18 Network visualization of keywords

Table 7 Top 1-,2-, and 3- word keywords used in the field

| # | 1-Word | Frequency | 2-Word | Frequency | 3-Word | Frequency |
|---|---|---|---|---|---|---|
| 1 | Optimization | 680 | Genetic Algorithms | 367 | Constraint-Handling Techniques | 74 |
| 2 | Algorithms | 360 | Constraint Handling | 184 | Particle Swarm Optimization (PSO) | 565 |
| 3 | Scheduling | 141 | Constrained Optimization | 1071 | Multiobjective Optimization | 467 |
| 4 | NSGA-II | 97 | Multiobjective Optimization | 1339 | Particle Swarm Optimization | 205 |
| 5 | Design | 96 | Evolutionary Algorithms | 1081 | Constrained multiobjective optimization | 68 |
| 6 | Algorithm | 59 | Differential evolution | 208 | Electric Load Dispatching | 66 |
| 7 | Reliability | 33 | Problem Solving | 239 | Multiobjective optimization problem | 222 |
| 8 | Investments | 33 | Multi objective | 307 | Differential evolution algorithms | 71 |
| 10 | Benchmarking | 113 | Decision making | 135 | Pareto optimal solutions | 121 |

| 12 | Costs | 57 | Pareto Principle | 281 | Constrained multiobjective optimizations | 213 |

## 6.8 Publication statistics by number of pages (pages count)

As of writing this paper, May of 2021, approximately 22395.8 pages of papers on constraint handling multiobjective population-based algorithms were published, with an average of 13.0435 pages per paper.

About 22.53931% of the articles possess between 10 and 15 pages; 12.17239% of the manuscripts are between 15 and 20 pages; 39.07979 % of the papers are between 5 and 10 pages; and 67.09377% of the manuscripts are between 5 and 20 pages. Figure 18 presents the distribution of the manuscripts based on page count.

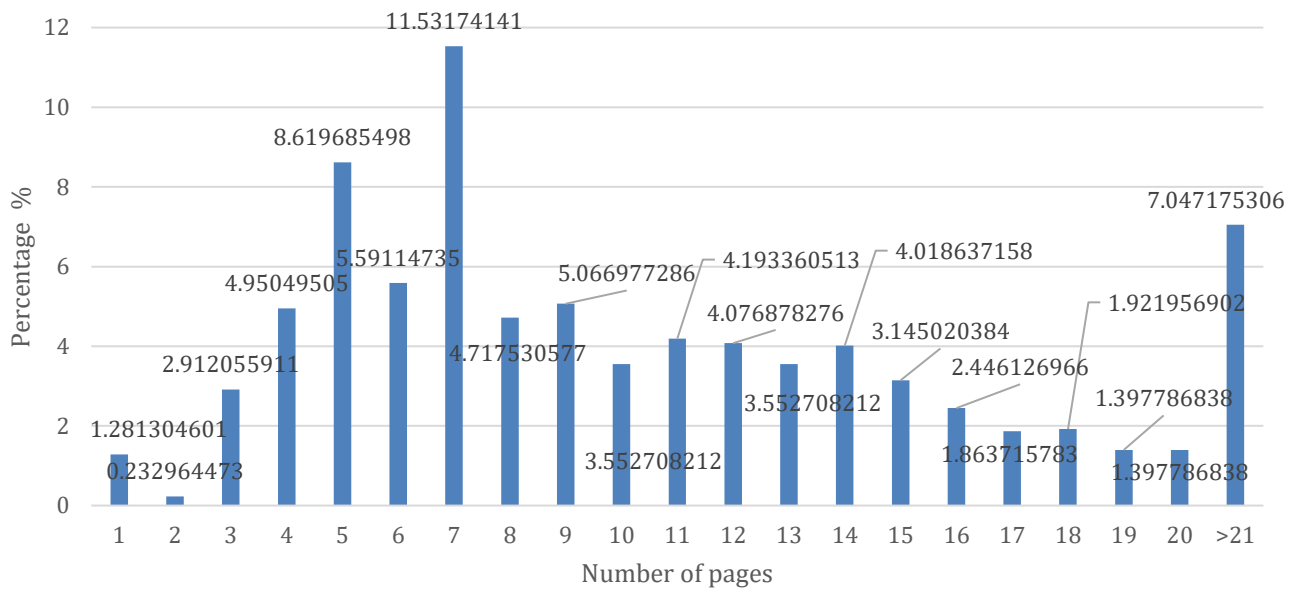

Figure 18 Distribution of documents based on page count

## 7 Summary and Future Research (RQ7)

The paper presents an analysis and overview of CHTs applied to multiobjective population-based algorithms. In the first part of the paper, the main idea of CHTs are defined, and the second part discusses a detailed scientometric analysis in the field. Some important technical points are extracted as follows:
- In the death penalty method, no information is used from infeasible points.
- Static penalty method is problem-dependent and may need several penalty parameters.
- Dynamic penalty method may converges to either an infeasible or feasible solution that is far from global optimum.
- The main disadvantages of the annealing penalty method is sensitivity to values of its factors.
- Setting parameters in adaptive penalty method is difficult and the method needs the definitions of additional parameters.
- In the self-adaptive penalty method, the additional parameters may affect the fitness function evaluations.
- The main difficulty in SGA is selecting the penalty factors for each sub-population.
- If the population is completely infeasible, choose the solutions that have a smaller overall constraint violation.
  - Retaining a proportion of infeasible solutions in the population may enhance the convergence and diversity of the algorithm.

- Two types of CHTs, namely repair methods and special genetic operators, focus only on the feasible space.
- Feasible solutions could be used to repair infeasible solutions (repairing population).
- According to the constraint dominance principle, the feasible solution is always prefered over the infeasible solution, which may cause loss of important information from infeasible individuals.
- Retaining a huge number of infeasible solutions may cause low convergence speed.
- Although special operators are known to be highly comparative CHT, their applicability is limited, which makes this technique difficult to run.
- Decoder is an interesting CHTs, but it involves a high computational cost and, thus, is now rarely used.
- Although the ensemble CHT has a competitive performerence, the method is parameter-dependent.
- Although the stochastic ranking method has been employed in several nature-inspired algorithms, it is not often used for the multiobjective version of the algorithms.
- Epsilon constraint method has been known as a powerful CHT, however, in some cases, premature convergence has been reported, while other works report that the method relies on gradient-based mutation.
- Using multiobjective concept as a CHT may require the gradient calculation.
- Recently, feasibility rules have been recognized as one of the most powerful CHTs, which are simple and flexible; however, one of the major disadvantages of this method is premature convergence since this technique favors feasible solutions.

As a future direction, the authors have identified the top 5 most-used keywords and research fields in the last three years (2019-2021) based on Scopus. Tables 9 and 10 show the mentioned keywords and research fields for this time period. It is obvious that multiobjective optimization, constraint optimization, and evolutionary algorithms are the most famous keywords in the last three years. It should be noted that CHTs for multiobjective optimization have not received much attention compared with single-objective optimization. It is suggested that researchers focus on such methods in future works. Also, BU technique, which is able to handle constraints directly possesses the potential to couple with a multiobjective evolutionary algorithm (MOEA) as well. Furthermore, it is suggested to focus on constraint handling techniques on many-objective optimization problems (with more than three objectives) as it is not received much attention. In addition, according to Tables 8 and 9, GA, DE, and PSO remain the top 3 algorithms, which are expected to be further explored in the future. Moreover, Engineering, Computer Science, and Mathematics have been the top 3 research fields in the last two years, and it is projected that research work will advance in these areas in the future. It is also recommended to review the applications of constrained multioobjective evolutionary algorithms in different sectors; including engineering design problems[248] [248][249][250], scheduling optimization problems [251]    [252] [253][253][254] , and resource optimization problems [255] [256].

*Table 8 Top 5 keywords in 2019 and 2020*

| # | Keywords (Scopus) | Frequency |
|---|---|---|
| 1 | Pareto principle | 65 |
| 2 | Genetic algorithms | 72 |
| 3 | Differential evolution | 45 |

| 4 | Particle swarm optimization (PSO) | 132 |
| 5 | Economic and social effects | 34 |
| 6 | Benchmarking | 32 |
| 7 | Decision making | 36 |
| 8 | Energy utilization | 28 |
| 9 | Scheduling | 39 |
| 10 | Pareto optimal solutions | 24 |

*Table 9 Top 5 research fields in 2019 and 2020*

| # | Research fields (Scopus) | (%) Contribution |
|---|---|---|
| 1 | Engineering | 24.3 |
| 2 | Computer Science | 31.8 |
| 3 | Mathematics | 17.2 |
| 4 | Energy | 5.9 |
| 5 | Decision Sciences | 4.1 |
| 6 | Materials Science | 4.1 |
| 7 | Business, Management and Accounting | 1 |
| 8 | Environmental Science | 2.1 |
| 9 | Physics and Astronomy | 2.8 |
| 10 | Earth and Planetary Sciences | 1.4 |

## 8 Discussion and Conclusion

Constraint population-based optimization involves the use of a population-based algorithm in combining with a CHT to solve a constraint optimization problem. This paper presents an analysis and evaluation of the CHTs on multiobjective optimization population-based algorithms, which support evolutionary and swarm intelligence algorithms. To the best of our knowledge, this study is the first analysis of relevant journals evaluated over the most relevant journals, keywords, authors, and articles in this field. All related papers, including research articles, reviews, book/book chapters, conference papers, etc., as of writing this paper, were extracted and analyzed. Publication statistics by year, journal, country, affiliation, author, number of pages, number of authors, and keywords are discussed in this paper as follows:

- According to WOS, 45824 citations have been received by the related papers, which is an average of 1992.35 citations per year and an average of 62.35 citations per item in WOS.
- Based on WOS, articles were the most popular document type, with a total of 522 articles (71.02%), and 2.60 authors per publication.
- Articles as the document type had the highest CPP2020 of 84.10, followed by proceedings papers with TP of 220 (29.93% of contributions and APP=2.13).
- Conference papers have the most contributions before 2010 followed by articles. However, since 2010, articles have the most contributions in the field.
- A total of 271 articles (36.87% of total), with 2.31 authors per publication (on average), were published in the category of Computer Science Artificial Intelligence, according to WOS.
- In total, 271 articles (36.87% of 735 articles) were published in the first category (Computer Science Artificial Intelligence) and a total of 83.39% were published in the first three categories: Engineering Electrical Electronic (23.26%) and Computer Science Theory Methods (23.26%).
- The highest CPP2020 of articles published in Computer Science Theory Methods is 190.011, which includes the paper "A fast and elitist multiobjective genetic algorithm: NSGA-II" by [16], and the highest APP for articles published in 'Energy fuels' is 2.97.
- "Computer Science Artificial Intelligence," "Engineering Electrical Electronic,"

"Computer Science Theory Methods," and "Computer Science Interdisciplinary Applications" were the top 4 productive WOS categories in the field.
- China, USA, and India were the top three active countries in the field according to WOS.
- Computer science, Engineering, and Mathematics have the most contributions with 936, 619, and 580 published articles, respectively. Pharmacology, Medicine, and Economics own the least contributions with 1, 4, and 6 published documents in the field, according to Scopus.
- Kalyanmoy Deb from "Michigan State University (USA)", Ray T. from "University of New South Wales (Australia)", and Carlos A. Coello Coello from "Cinvestav-IPN (Mexico)" are the top 3 authors in the field with 38, 32, and 28 publications (indexed by Scopus), respectively. WANG Y from "City University Hong Kong (Hong Kong)", Carlos A. Coello Coello from "Cinvestav-IPN (Mexico)", and Ray T. from "University of New South Wales (Australia)" are the top 3 authors in the area with 21, 20, and 17 documents (indexed by WOS), respectively.
- Almost 0.5241% of authors own more than 10 documents; 1.6307% possess between 5 and 10 documents; 4.5428% have between 3 and 5 papers; 11.7647% of authors own 2 papers; and 81.5377% of authors possess 1 document.
- Approximately 22.53931% of the articles possess between 10 and 15 pages; 12.17239% of the manuscripts are between 15 and 20 pages; 39.07979% of the papers are between 5 and 10 pages; and 67.09377% of the manuscripts are between 5 and 20 pages.

**Author Contributions:** "Conceptualization, Iman Rahimi. and Amir H. Gandomi.; methodology, Iman Rahimi; software, Iman Rahimi; validation, Amir H. Gandomi, Fang Chen. and Efrén Mezura-Montes; formal analysis, Iman Rahimi; data curation, Iman Rahimi; writing—original draft preparation, Iman Rahimi—review and editing, Amir H. Gandomi, Fang Chen. and Efrén Mezura-Montes ;supervision, Amir H. Gandomi, Fang Chen;

**Funding:** "This research received no external funding".

**Conflicts of Interest:** "The authors declare no conflict of interest."